\documentclass[11pt]{article}
  
\newcommand{\di}{\displaystyle}
\newcommand{\iy}{\infty}
\newcommand{\lt}{\left}
\newcommand{\me}{\medskip}

\newcommand{\ri}{\rightarrow}
\newcommand{\rt}{\right}
\newcommand{\sm}{\smallskip}

\newcommand{\wi}{\widetilde}
\newcommand{\wit}{\widehat}

\newcommand{\ex}{\exists\ }
\newcommand{\fo}{\forall\ }

\newcommand{\Id}{\mathrm{Id}}

\newcommand{\st}{\,:\,}

\newcommand{\un}{\text{\bf 1}}

\newcommand{\bq}{\begin{eqnarray*}}
\newcommand{\bqn}[1]{\begin{eqnarray}\label{#1}}
\newcommand{\eq}{\end{eqnarray*}}
\newcommand{\eqn}{\end{eqnarray}}
\newcommand{\wwtbp}{\par\hfill $\blacksquare$\par\me\noindent}
\newcommand{\df}{:=}

\usepackage[latin1]{inputenc}

\usepackage{amsfonts}      
\usepackage{amsmath}     
\usepackage{amssymb} 

\newtheorem{theo}{Theorem}[section]
\newtheorem{prop}[theo]{Proposition}               
\newtheorem{lem}[theo]{Lemma}           
\newtheorem{cor}[theo]{Corollary}

\newtheorem{rem}[theo]{Remark}             
\newtheorem{defi}[theo]{Definition}          
\newtheorem{nota}[theo]{Notation}
\numberwithin{equation}{section}

\newcommand{\EE}{\mathbb{E}}

\newcommand{\NN}{\mathbb{N}} 
\newcommand{\PP}{\mathbb{P}} 
 
\newcommand{\RR}{\mathbb{R}}

\newcommand{\Aa}{ {\cal A }} 
\newcommand{\Ba}{ {\cal B }} 
\newcommand{\Ca}{ {\cal C }} 
\newcommand{\Da}{ {\cal D }} 
\newcommand{\La}{ {\cal L }}

\newcommand{\Ea}{ {\cal E }} 
\newcommand{\Sa}{ {\cal S }}

\newcommand{\Fa}{ {\cal F }} 
\newcommand{\Ga}{ {\cal G }} 
\newcommand{\Qa}{ {\cal Q }} 
 
\newcommand{\Xa}{ {\cal X }} 
\newcommand{\Ma}{ {\cal M }} 
\newcommand{\Ta}{ {\cal T}}

\newcommand{\Pa}{ {\cal P }} 
\newcommand{\Za}{ {\cal Z }}

\newcommand{\Bb}{ {\bf B }} 
\newcommand{\Cb}{ {\bf C }}

\newcommand{\bb}{ {\bf b }} 
\newcommand{\ab}{ {\bf a }} 
\newcommand{\cb}{ {\bf c }} 
\newcommand{\pb}{ {\bf p }} 
\newcommand{\ib}{ {\bf i }} 

\newcommand{\point}{\mbox{\LARGE .}} 
 
\newcommand{\cqfd}{\hfill\blbx \\} 
\def\blbx{\hbox{\vrule height 5pt width 5pt depth 0pt}\medskip} 
 
\def \PP{\mathbb{P}} 
\def \RR{\mathbb{R}} 
\def \EE{\mathbb{E}}

\begin{document}

\title{The convergence to equilibrium of neutral genetic models} 
\author{P. Del Moral\footnote{CNRS UMR 6621, 
Universit\'e de Nice, Laboratoire de Math\'ematiques J. Dieudonn\'e,  
Parc Valrose, 
06108 Nice, France; and IRISA / INRIA - Campus universitaire de Beaulieu - 35042 Rennes, France. }, 
L. Miclo\footnote{Laboratoire d'Analyse, Topologie, Probabilit\'es UMR 6632, Universit\'e de Provence, Technop\^ole Ch\^ateau-Gombert, 39, rue F. Joliot Curie, 13453 Marseille Cedex 13
},
F. Patras\footnote{CNRS UMR 6621, 
Universit\'e de Nice, Laboratoire de Math\'ematiques J. Dieudonn\'e,  
Parc Valrose, 
06108 Nice Cedex 2, 
France. The second author was supported by the ANR grant AHBE 05-42234},
S. Rubenthaler\footnote{CNRS UMR 6621, 
Universit\'e de Nice, Laboratoire de Math\'ematiques J. Dieudonn\'e,  
Parc Valrose, 
06108 Nice Cedex 2, 
France}}

\maketitle 
 
\begin{abstract}
This article is concerned with the long time behavior of neutral
genetic population models, with fixed population size. We design 
an explicit, finite, exact, genealogical tree based representation of stationary
populations that holds both for finite and infinite types (or alleles) models.
We then analyze the
decays to the equilibrium of finite populations in terms of the convergence
to stationarity of their first common ancestor. We estimate the Lyapunov exponent of the distribution flows with respect to the total variation norm. We give bounds on these exponents only depending on the stability with respect to mutation of a single individual; they are inversely proportional to 
the population size parameter.
\

\


\end{abstract}

\section{Introduction}

Stochastic models for population dynamics provide a mathematical framework for the analysis of genetic variations in biological populations that evolve under the influence of evolutionary type forces such as selections and mutations. 
The common models, such as the (generalized) Wright-Fisher models, allow for many variations in the fundamental assumptions. For example, one may consider  a finite population, or a sample out of an infinite population; there might be a finite or infinite number of types (or alleles) of individuals; the genetics involved may refer to monoecious (i.e. where there is only one sex) or dioecious phenomena, and so on.

Dealing with finite number of individuals without sampling in an infinite population is usually difficult: although the problems are easily formulated, computations become soon very involved.
Here, we are interested in finite population models where the evolution is driven by a selection/mutation process and derive explicit formulas for the stationary state of the population or the decay to the equilibrium. We do not put any restriction on the number of types or alleles, that may be finite or infinite. The main restriction to a full generality of the model is neutrality of the selection process: that is, we assume that all individuals have the same reproduction rate (we refer e.g. to the monograph of  M. Kimura~\cite{kimura} for a detailed account on the neutral theory of molecular evolution).

At the molecular level, the evolutionary models we consider
correspond therefore to situations where
the genetic drifts that govern the dynamics are the mutations combined with a uniform reproduction rate. Genealogical tree evolution models arise then in a natural way, when considering the complete past history of the individuals. These path space models can be described  forwards with respect to the time parameter. Conversely, we can also trace back in 
time the complete ancestral line of all the individuals. This backward view of the ancestral structures
is then interpreted as a
stochastic coalescent process. 

The important questions that arise then are the precise description of the asymptotic behavior of  evolution processes, and the corresponding time to equilibrium analysis. They are at the core of the modern development of mathematical biology. Let us also point out that, apart from their importance in biology, they
are also related to the convergence of a class of genealogical tree algorithms used in advanced stochastic 
engineering, and in Bayesian statistics. For example, a full understanding 
of the long time behavior of genetic type branching models 
is essential in the designing of
genealogical particle filters and smoothers, as well as
for the tuning of genetic algorithms for solving global
optimization problems. These features of evolutionary models are another strong motivation for the present work and, although we emphasize mainly the applications to genetics, the interested reader should keep in mind these other application areas (for further informations and references on the subject, see for instance the pair of recent books~\cite{fk,doucet}). 

 During the last three decades, 
many efforts have been made to tackle these
questions.  Although several natural Markov chain models of genetic processes have been developed, their combinatorial complexity makes both the intuitive, and the rigorous understanding of the long time behavior of evolution mechanisms difficult. Several reduced models have been developed to obtain rigorous, but partial theoretical
results. These simplified models are usually based on 
large population approximations, and appropriate rescaling
of the time parameter in units of the total population size.

Let us review very briefly three main directions of research in the field.  The first one studies the diffusion limit model arising,
 through appropriate time
rescaling techniques, when considering gene type
fractions decompositions on large population sizes models (see for instance the seminal book
of S. N. Ethier and T. Kurtz~\cite{kurtz}, and the references therein).
A second idea is to consider the backward ancestral lineage in infinitely many allele models. This genealogical process is expressed in terms of rescaled Poisson coalescence epochs, and superimposed 
 Poisson mutation events. Running back in time this population model, the coalescent of Kingman describes the ancestral genealogy of the
 individuals in terms of a simplified binary ancestral lineage tree.
 The Ewens' sampling formula describes the limiting distribution of
 the type spectrum in finite samples of the infinitely many allele models (see for instance the Saint Flour's lecture notes of S. Tavar\'e~\cite{tavare}, and references therein). A third, more recent, idea is essentially based on the mean field interpretation
of genetic models.  In this context, the occupation measures
of simple genetic populations converge, as the size of the systems
tends to infinity, to a non linear
 Feynman-Kac semigroup in the distribution space. For a rather detailed account
 of these mean field limits, we refer the reader to the first author monograph~\cite{fk}. In this interpretation, and in the case of reversible mutations, the long time behavior
 of an infinite population model corresponds to the ground state of Schr\"odinger type operators. 
 
In the present article we depart from these techniques and tackle directly the combinatorial complexity of finite population dynamics. This leads to a fine asymptotic analysis of a general class of neutral genetic models on arbitrary state spaces, with fixed
population size. In contrast to the existing literature 
on the asymptotic behavior of neutral genetic models, our approach is
not based on some had-oc rescaling of the time parameter, and it 
applies to evolutionary models with fixed population size. 

Firstly, we provide an 
explicit functional representation of their invariant measures in terms of planar genealogical trees. We then use this representation to analyze  the decays to stationary populations in terms of the convergence to the equilibrium of the first common ancestor and show, for example, that the
Lyapunov exponent of the genetic distribution flow is inversely
proportional to the population size of the model.
This estimate is of course sharper than 
the one we would obtain using crude minorization techniques of the 
transitions probabilities of the genetic model, such as those presented in~\cite{gao}. We refer to Section \ref{asbehav}, Theorems \ref{theo1} and \ref{theo2} for a precise statement of the main results of the article.

Unfortunately, these rather natural genealogical techniques 
are restricted to neutral genetic population models, but it leads to conjecture that the same type of decays holds true for more general models. Besides, some of the ideas and techniques developed in the article can certainly be adapted to cover inhomogeneous selection rates -this is a point we leave deliberately out of the scope of the present work.

We finally mention that the genealogical invariant measure derived in the present article belongs to the same class of tree based measures as the ones studied in~\cite{dpr} to derive coalescent tree based functional representations of genetic type particle models. The symmetry properties of the invariant measure
can be combined with the combinatorial analysis presented in~\cite{dpr} to simplify notably the formula for the distributions of stationary populations. We briefly present these constructions in the last section.

To the best of our
knowledge, these genealogical tree based representations
of stationary populations, as well as of the corresponding convergence
decays to the equilibrium, are the first of this kind, for this class of 
evolutionary models.

\section{Neutral Genetic Models}
\subsection{Conventions}
Let us introduce some notation. We denote respectively
by $\mathcal{M}(E)$,   $\mathcal{P}(E)$, and $\mathcal{B}(E)$, the set
of all finite signed measures on some measurable space
$(E,\mathcal{E})$,  the convex subset
of all probability measures,
and
 the Banach
space of all bounded and measurable functions $f$ on $E$, equipped with
the uniform norm $\Vert f\Vert=\sup_{x\in E}{|f(x)|}$. The space $E$ will stand for the set of types or alleles of the genetic model. 

We let $
\mu(f)=\int~\mu(dx)~f(x)
$, be  the
 Lebesgue integral of a function  $f\in\mathcal{B}(E)$, with respect to
a measure $\mu\in \mathcal{M}(E)$, and we equip $\Ma(E)$ with the total variation norm
$\|\mu\|_{\rm tv}=\sup_{f\in\Ba(E):\Vert f\Vert\leq 1}| \mu(f)|$. 

We recall that a Markov transition $M$ from $E$
into itself, is an integral probability operator  such that
 the functions
$$
M(f)~:~x\in E\mapsto M(f)(x)=
\int_{E}~M(x,dy)~f(y)~\in \RR
$$
are $\Ea$-measurable and bounded,  for any $f\in\mathcal{B}(E)$. It generates a dual operator $\mu\mapsto \mu M$
from $\mathcal{M}(E)$ into itself
defined by $(\mu M)(f):=\mu(M(f))$. For a pair of Markov transitions  $M_1$, and $M_2$,
we
denote by $M_1M_2$ the composition integral operator from
$E$ into itself, defined for any $f\in\Ba_b(E)$ by
$(M_1M_2)(f):=M_1(M_2(f))$. The tensor power $M^{\otimes N}$
represents
the bounded integral operator  on $E^N$, defined  for any $F\in\mathcal{B}(E^N)$ by
 $$
 M^{\otimes N}(F)(x^N,\ldots,x^N)=\int_{E^N}~\left[M(x^1,dy^1)~\ldots~M(x^N,dy^N)\right]~F(y^1,\ldots,y^N)
 $$

 We  let
$N$ be a fixed integer parameter, and we set $[N]=\{1,\ldots,N\}$. The integer $N$ will stand for the number of individuals in the population.
We denote by $m(x)=\frac{1}{N}\sum
_{i=1}^{N}\delta_{x^{i}}$ the empirical measure associated with an $N$-uple $x=(x^i)_{1\leq i\leq N}\in E^N$.
The notation $\delta_x$ stands in general for the Dirac measure 
at the point $x$. 

We let $\Aa=[N]^{[N]}
$, be the set of mappings from $[N]$ into itself. We associate with a
mapping $a\in\Aa$, the Markov transition $D_{a}$ on $E^N$ defined by
$D_{a}(x,dy)=\delta_{x^a}(dy)$, with $x^a=(x^{a(i)})_{i\in[N]}$. The transition $\Da_a$ can be interpreted as a coalescent, or a selection type transition on the set $E^N$. In this interpretation, the 
population of individuals $x^a=(x^{a(i)})_{i\in[N]}$ results
from a selection in $x=(x^{i})_{i\in[N]}$ of the individuals with
labels $(a(i))_{i\in[N]}$. Also notice that the $N$-tensor product of an empirical measure
can be expressed in terms of these selection type transitions :
$$
m(x)^{\otimes N}(dy)=\frac{1}{|\Aa|}\sum_{a\in\Aa}D_{a}(x,dy)=:\Da(x,dy)$$

We consider now a Markov transition  $M$, and a probability measure $\eta$ on $E$. The Markov transition $M$ will represent the mutation transition process, whereas $\eta$ will stand for the initial distribution of the population.
An $N$-neutral genetic model, with these parameters, can be represented by a pair of 
$E^{N}$-valued Markov chains $(\xi_n)_{n\in\NN}=((\xi_{n}^{i})_{i\in [N]})_{n\in\NN}$ and 
$(\widehat{\xi}_n)_{n\in\NN}=((\widehat{\xi}_{n}^{i})_{i\in [N]})_{n\in\NN}$, together with transitions
\begin{equation}\label{defgen}
\xi_n\stackrel{ selection}{-\!\!\!\!-\!\!\!\!-\!\!\!\!-\!\!\!\!-\!\!\!\!-\!\!\!\!-\!\!\!\!-\!\!\!\!\longrightarrow}\widehat{\xi}_n
\stackrel{ mutation}{-\!\!\!\!-\!\!\!\!-\!\!\!\!-\!\!\!\!-\!\!\!\!-\!\!\!\!-\!\!\!\!-\!\!\!\!\longrightarrow}\xi_{n+1}
\end{equation}
The initial
configuration $\xi_{0}$ consists of $N$ independent
and identically distributed random variables with common distribution $\eta$.  

The selection transition $\xi_n\leadsto\widehat{\xi}_n$ corresponds to a simple Wright-Fisher selection model: it consists in sampling $N$ conditionally independent random variables $(\widehat{\xi}_n^i)_{i\in[N]}$, with common distribution $m(\xi_n)$. Equivalently, 
the selected population $\widehat{\xi}_n$ is a random variable with distribution $\Da(\xi_n,dx)$; in other words, we chose randomly a mapping $A_n$ in $\Aa$, and we set 
$\widehat{\xi}_n=\xi_n^{A_n}\left(=(\xi_n^{A_n(i)})_{i\in[N]}\right)
$. In the literature on branching processes, this random selection mapping $A_n$ 
is often expressed in terms of 
a multinomial branching rules. Assume,  for instance that all the $\xi_n^i$ are different, if we let $L_n=(L^i_n)_{i\in[N]}$ be the number of offsprings
of the individuals $\xi_n$
$$
\forall i\in[N]\qquad L^i_n:=\left|\left\{j\in [N]~:~\widehat{\xi}_n^j=\xi^i_n\right\}\right|
$$
then, we find that $L_n$ has a symmetric multinomial distribution
$$
\PP\left(L_n^1=l^1,\ldots,L_n^N=l^N~|~\xi_n\right)=\frac{N!}{l^1!\ldots l^N!}~\frac{1}{N^{N}}
$$
for any $l=(l^i)_{i\in[N]}\in\NN^N$, with $\sum_{i\in[N]}l^i=N$. 

During the mutation, each of the selected particles $\widehat{\xi}_n^i$ evolves to a new location $\xi_{n+1}^i$, randomly chosen with distribution $\delta_{\widehat{\xi}_n^i}M$, with $i\in[N]$. Equivalently,
the population ${\xi}_{n+1}$ after mutation is a random variable with distribution $\Ma(\widehat{\xi}_n,dx)$, with $\Ma:=M^{\otimes N}$.

The distribution laws of the populations before and after 
the selection step are defined by the values, for any function $F\in\Ba_b(E^N)$, of
$$
\Gamma_{\eta,n}(F):=\EE\left(F(\xi_{n})\right)
\quad\mbox{\rm and}\quad
\widehat{\Gamma}_{\eta,n}(F):=\EE\left(F(\widehat{\xi}_{n})\right)
$$
Notice that $\Gamma_{\eta,0}=\eta^{\otimes N}$, and $\widehat{\Gamma}_{\eta,0}=\Gamma_{\eta,0}\Da$. By construction, we also have that
$
\widehat{\Gamma}_{\eta,n}=\Gamma_{\eta,n}\Da$, and
$
\Gamma_{\eta,n+1}=\widehat{\Gamma}_{\eta,n}\Ma
$.
This clearly gives the dynamical structure of the pair of distributions
$$
\Gamma_{\eta,n+1}=\Gamma_{\eta,n}\left(\Da\Ma\right)
\quad\mbox{\rm and}\quad
\widehat{\Gamma}_{\eta,n+1}=\widehat{\Gamma}_{\eta,n}\left(\Ma\Da\right)
$$
In particular, we have:
$$
{\Gamma}_{\eta,n}(F)={\Gamma}_{\eta,0}\left(\Da \Ma\right)^n(F)=
\EE\left(\Gamma_{\eta,0}D_{A_0}\Ma D_{A_1}\ldots D_{A_{n-1}}\Ma (F)\right)
$$
and
$$
\widehat{\Gamma}_{\eta,n}(F)=\widehat{\Gamma}_{\eta,0}\left(\Ma \Da\right)^n(F)=
\EE\left(\Gamma_{\eta,0}D_{A_0}\Ma D_{A_1}\ldots \Ma D_{A_{n}}(F)\right)
$$
with a sequence of independent random
variables $(A_p)_{0\leq n}$, uniformly chosen in the set $\Aa$.

\subsection{Genealogic trees representations}
In neutral reproduction models, the pair of mutation-selection processes can be separated in order to describe the 
mutation scenarios and the neutral selection reproductions on two different levels. Traditionally, these neutral genetic models
are sampled in the following way. The genealogical
structure of the individuals is first modelized. Then, the genetic mutations are superimposed 
on the sampled genealogy. 

In the present section, we design a representation of this pair of genetic processes in terms of random trees.
The neutral selections  are encoded by random
mappings; their composition gives rise to random coalescent trees. 
The representation of the process is complete when mutation scenarios are taken into account. In contrast to traditional mutation-selection decoupling models, non necessarily neutral selection processes could also be described 
forward 
in terms of these random trees. 
In this situation, the random selection type mappings would depend on the
configuration of the genes.

We let $X_n^{\ib_n}$, with $\ib_n=(i_0,\ldots,i_n)\in[N]^{n+1}$, and $n\geq 0$,  be the collection of $E$-valued random variables defined inductively as follows: at rank $n=0$, $X_0^{\ib_0}$, with $\ib_0=i_0\in[N]$, represents a collection of $N$ independent, and identically distributed random variables with common distribution $\eta$. 
Given the random variable $X_{n-1}^{\ib_{n-1}}$, 
for some multi index $\ib_{n-1}\in[N]^{n}$, the sequence of random variables 
$X_n^{\ib_n}$, with $\ib_n=(\ib_{n-1},i_n)$, and $i_n\in[N]$ consists in $N$ conditionally independent random variables, with common distribution $\delta_{X_{n-1}^{\ib_{n-1}}}M
$.

Following the convention of~\cite{harris1}, the collection of
random variables $X_{p}^{\ib_{p}}$, with $\ib_p\in[N]^{p+1}$, and $0\leq p\leq n$, 
can be associated to the vertices of a planar forest of height $n$. 
The
sequence of integers $\ib_n$ represents the complete genealogy of an individual at level $n$. For instance, an individual $X_{2}^{\ib_{2}}$ in the second generation is 
associated with the triplet $\ib_2=(i_0,i_1,i_{2})$
of integers that indicates that he his the $i_2$-th child of the $i_{1}$-th child of the $i_{0}$-th root ancestor individual. Running back in time, we can trace back the complete ancestral line of a given current individual
$X_n^{\ib_n}$ at the $n$-th generation
$$
X_0^{\ib_0}\longleftarrow X_1^{\ib_1}\longleftarrow
\ldots\longleftarrow
X_{n-1}^{\ib_{n-1}}\longleftarrow
X_n^{\ib_n}
$$
The neutral selection, and the mutation transition introduced in (\ref{defgen}) have a natural interpretation in terms of these random forests. Roughly speaking, the random forests introduced above represent all 
the transitions of  a selection-mutation genetic algorithm.  

To be more precise, it is convenient to introduce some additional algebraic structures. We equip  $\Aa=[N]^{[N]}$ with the unital semigroup structure
associated with the 
composition operation
$ab:=a\circ b$, and the identity element $Id$.
We equip
the set $\Aa$ with the partial order relation defined for any
pair of mappings $a,c\in \Aa $ by the following 
formula
$$
a\leq c\Longleftrightarrow
\exists b\in\Aa\quad a=b c
$$
For any collection of mappings $(a_p)_{0\leq p\leq n}$, we notice that
$$
a_{0,n}\leq a_{1,n}\leq\ldots\leq a_{n-1,n}\leq a_{n}\leq Id 
$$
with the composition semigroup $a_{p,n}=a_{p}a_{p+1,n}$, $0\leq p< n$; and the convention $a_{n,n}=a_n$.
For any weakly decreasing sequence of mappings  $(b_0,\ldots,b_n)\in\Aa^{n+1}$ (that is, any sequence s.t. $b_i\geq b_{i+1}$, we set
$$
\Xa^{(b_n,\ldots,b_0)}_n:=\left(X_n^{b_{n}(i),\ldots,b_1(i),b_0(i)}\right)_{1\leq i\leq N}
$$
In this notation, it is now easy to prove  that the
neutral genetic model can be seen as a particular
way to explore the random forest introduced above 
$$
\Xa_0^{Id}\stackrel{ selection}{-\!\!\!\!-\!\!\!\!-\!\!\!\!-\!\!\!\!-\!\!\!\!-\!\!\!\!-\!\!\!\!-\!\!\!\!\longrightarrow}\Xa_0^{A_0}\stackrel{ mutation}{-\!\!\!\!-\!\!\!\!-\!\!\!\!-\!\!\!\!-\!\!\!\!-\!\!\!\!-\!\!\!\!-\!\!\!\!\longrightarrow}
\Xa_{1}^{A_0,Id}
\stackrel{ selection}{-\!\!\!\!-\!\!\!\!-\!\!\!\!-\!\!\!\!-\!\!\!\!-\!\!\!\!-\!\!\!\!-\!\!\!\!\longrightarrow}
\Xa_1^{A_0A_1,A_1}\stackrel{ mutation}{-\!\!\!\!-\!\!\!\!-\!\!\!\!-\!\!\!\!-\!\!\!\!-\!\!\!\!-\!\!\!\!-\!\!\!\!\longrightarrow}
\Xa_2^{A_0A_1,A_1,Id}\longrightarrow
\ldots
$$
with a sequence of independent random
variables $(A_p)_{0\leq n}$, randomly chosen in the set $\Aa$.
More generally, the distribution laws of the neutral genetic model
are given by the formulae
$$
\widehat{\Gamma}_{\eta,n}(F)=\EE_{\eta}\left(F\left(\Xa_n^{A_{0,n},A_{1,n},\ldots,A_{n-1,n},A_n}\right)\right)
$$
with 
the semigroup
$A_{p,n}:=A_pA_{p+1,n}$, with $0\leq p< n$,
and the convention $A_{n,n}=Id$. In much the same way, we also have that
$$
\Gamma_{\eta,n+1}(F)=\widehat{\Gamma}_{\eta,n}\Ma(F)=\EE_{\eta}\left(F\left(\Xa_{n+1}^{A_{0,n},A_{1,n},\ldots,A_{n-1,n},A_n,Id}\right)\right)
$$
\subsection{Asymptotic behavior}\label{asbehav}
  The permutation mappings $\sigma\in\Ga_{N}$ are the largest elements of $\Aa$, while
the smallest elements of $\Aa$ are given by the constants elementary mappings  $e_i(j)=i$, for any $j\in[N]$, with $1\leq i\leq N$. 
The set of these lower bounds is denoted by
$$
\partial \Aa=\left\{a\in\Aa~:~|a|=1\right\}=\{e_i~:~i\in[N]\},
$$
where $|a|$ stands for the number of elements in the image of $a$.
We also let $(B_p)_{0\leq p\leq n}$ be the weakly decreasing Markov chain on $\Aa$ defined by
\begin{equation}\label{defb}
\forall n\geq 1\qquad B_n:=A_nB_{n-1}
\end{equation}
with the initial condition $B_0=A_0$, where $(A_p)_{0\leq p\leq n}$ stands for a sequence of independent and uniformly distributed random variables on the set of mappings $\Aa$. We let $T$ be the first time the chain $B_n$ enters in the set $\partial \Aa$.
\begin{equation}\label{defT}
T:=\inf{\left\{n\geq 0~:~B_n\in \partial \Aa\right\}}
\end{equation}
In terms of genealogy, the Markov chain (\ref{defb})
represent the ancestral branching process from the present generation at time $0$ back into the past. In this interpretation, 
the mapping $A_n$ in formula (\ref{defb}) represents the way the
individuals choose their parents in the previous ancestral generation. 
The  range of $A_n$ represents successful parents with direct descendants, whereas the range of $B_n$ represents successful ancestors with descendants in all the generations till time $0$. The random variable $T$ represents the time to most recent common ancestor of an initial population with $|B_0|$ individuals.
 
We are now in position to state the two main results of this article. 
\begin{theo}\label{theo1}
The Markov chain $B_n$ is absorbed by the boundary $\partial\Aa$
in finite time, and $B_T$ is uniformly
distributed in $\partial\Aa$.  In addition, for  any time horizon $n\geq 0$, we have
$$
\PP(T>n)\leq
K \left(\frac{n}{N}\vee 1\right)\exp\left( -\left(\frac{n}{N}-1\right)_+ \right)\quad\mbox{\rm with}
\quad 
K:=e \prod_{l\in\NN^*} \left(1-\frac{2}{(l+1)(l+2)}\right)^{-1}
\quad 
$$
\end{theo}
The Theorem follows from (\ref{majtemps}).

\ \

We further assume that the Markov transition $M$ has an invariant probability measure
$\mu$ on $E$, and we denote by $\widehat{\Gamma}_{\mu}$ the
probability measure on $E^N$ defined by
\begin{equation}\label{measdef}
\forall F\in\Ba_b(E^N)\qquad
\widehat{\Gamma}_{\mu}(F):=\EE_{\mu}\left(F\left(\Xa_T^{B_T,\ldots,B_1,B_0}\right)\right)
\end{equation}
\begin{theo}\label{theo2}
The measure $\widehat{\Gamma}_{\mu}$ is an invariant measure of the $N$-neutral genetic model $\widehat{\xi}_n$. In addition, if the mutation transition $M$ satisfies the following regularity condition :

{\rm\bf (H)} {There exists some  $\lambda>0$, and some finite constant $\delta<\infty$ such that for any $n\geq 0$
\begin{equation}\label{dob}
\beta(M^n):=\sup_{(x,y)\in E^2}\|M^n(x,\point)-M^n(y,\point)\|_{\rm tv}\leq \delta~e^{-\lambda n}
\end{equation}}
then, there exists $K'\geq 0$ (a universal constant) such that  for  $n\geq N+\frac{1}{\lambda}$, 
we have the estimate
\begin{eqnarray*}
\left\|\widehat{\Gamma}_{\eta,n}-\widehat{\Gamma}_{\mu}\right\|_{\rm tv}&\leq& \delta K'\frac{n}{N}
\exp\left(-\frac{n}{N+\frac{1}{\lambda}}\right)
% +2K\frac{n}{N}\exp\left(-\left(\frac{n}{N}-1\right)_+\right)
\end{eqnarray*}
\end{theo}
This Theorem follows from Theorem \ref{theo1}, Equation (\ref{inter}) and Corollary \ref{majexp}. Notice that the regularity condition is a weak condition, most often satisfied in selection-mutation genetic models -degenerate cases should be avoided such as, for example, in the finite state case, the splitting of the matrix representation of $M$ into a direct sum, or the existence of several fixed points of the Markov transition. The regularity hypothesis is satisfied for example if $M(x,dy)=\alpha \delta_x(dy)+(1-\alpha )\mu (dy)$ with $0<\alpha <1$ and $\mu$ a Lebesgue-type density on $E$.

We end this section with some consequence of these two theorems.
Firstly, we notice that the first assertion of the theorem can alternatively be expressed in terms of the measure $\Gamma_{\mu}$ defined by
$$
\forall F\in\Ba_b(E^N)\qquad
\Gamma_{\mu}(F):=\EE_{\mu}\left(F\left(\Xa_T^{B_T,\ldots,B_1,B_0,Id}\right)\right)
$$
More precisely we notice that, by construction,
$\Gamma_{\mu}(F)=\widehat{\Gamma}_{\mu}(\Ma F)$. This readily implies that $\Gamma_{\mu}$ is an invariant measure of the neutral genetic model $\xi_n$. Indeed, we have
$$
\Gamma_{\mu}(\Da\Ma(F))=\widehat{\Gamma}_{\mu}(\Ma\Da(\Ma F))=
\widehat{\Gamma}_{\mu}\Ma (F)=\Gamma_{\mu}(F)
$$
The reverse assertion follows the same line of arguments.

Notice also that the second assertion of the theorem can be used to estimate the Lyapunov exponent of the  distribution semigroup of the neutral genetic particle model; that is we have that
$$
\liminf_{n\rightarrow\infty}-\frac{1}{n}\log{\left\|\widehat{\Gamma}_{\eta,n}-\widehat{\Gamma}_{\mu}\right\|_{\rm tv}}\geq \frac{\lambda}{\lambda N+1} \wedge \frac{1}{N}
$$

\section{Random mappings and coalescent tree based measures}

Recall that we denote
 by $|a|$ we denote the cardinality of the set $a([N])$. 
Notice that $a\leq c\Rightarrow |a|\leq |c|$, and the cardinality of the set 
$$
\Aa_{c}:=\{a\in\Aa~:~a\leq c\}
$$
coincides with the number  of ways $N^{|c|}$ of mappings the range
$c([N])$ of the mapping $c$ into the set $[N]$. 
Also observe that $\left|\partial \Aa\right|=N$, and $\Aa_{e_i}=\partial\Aa$, for any $i\in[N]$, where $|\delta A|$ stands for the cardinality of $\delta A$. In this notation, the transitions probabilities
of the chain $B_n$ introduced in (\ref{defb}) are given by

\begin{equation}\label{condexp}
\PP(B_n=b~|~B_{n-1}=a)=\frac{1}{N^{|a|}}~1_{\Aa_a}(b)
\end{equation}

It is also readily checked that the uniform distribution $\nu_0(a):=\frac{1}{|\Aa|}$
on $\Aa$ is such that $\nu_0(a)=K(Id,a)$, and we also have that
$K(e_i,b')=\frac{1}{N}~1_{\partial\Aa}(b')$, where $K$ stands for the Markov transition on $\Aa$ defined by 
$$K(\phi ):a\in\Aa\longmapsto K(\phi )(a):=\frac{1}{|\Aa |}\sum\limits_{b\in\Aa}\phi (ba)$$
for any function $\phi$ on $\Aa$.
 This clearly implies that the uniform measure $\nu_1(a):=\frac{1}{N}1_{\partial A}(a)
$ on the set $\partial A$ is an invariant measure of $K$. That is we have that $\nu_1=\nu_1K$.

 For any weakly decreasing sequence of mappings $\bb=(b_0,\ldots,b_n)\in\Aa^{n+1}$, with a finite length $l(\bb)=n$,  we write 
 $$
  |\bb|:=\sum_{0\leq p\leq n}|b_p|$$
  For any $\cb=(c_0,\ldots,c_n)\in\Aa^{n+1}$, any mapping $a\in\Aa$, and any $b\geq c_0$, we also write 
 $$
 (b,\cb)=(b,c_0,\ldots,c_n)
\qquad
\cb a:=(c_0a,\ldots,c_na)\quad \mbox{\rm and}\quad
\cb\star a:= (c_0a,\ldots,c_{t_n}a)
$$
where $t_n=\inf{\{0\leq p\leq n~:~c_pa\in\partial \Aa\}}$; with the convention
$t_n=n$, if $c_pa\not\in\partial \Aa$ for any $0\leq p\leq n$. 

The set of all weakly decreasing excursions from $\Aa$ into $\partial\Aa$ is given by
$$
\Ca:=\cup_{n\geq 0}\Ca_n
$$
with  the sets $\Ca_n$ of all weakly decreasing
 excursions
$\cb$,
with length $l(\cb)=n\geq 0$, and defined by
$$
\Ca_n=\left\{
\cb=(c_0,\ldots,c_n)\in\left(\Aa-\partial\Aa\right)^{n}\times\partial\Aa~:~\forall 0\leq p<n~~ c_p\geq c_{p+1}\right\}
$$
We use the convention $\Ca_0=\partial\Aa$, for $n=0$. We associate with the random excursion
$$
\Bb=(B_0,B_1,\ldots,B_T)\in\Ca
$$ 
the pair of random excursions
$$
\Bb'=(Id,\Bb)\in\Ca
\quad\mbox{\rm and}\quad
\Bb' \star A\in\Ca
,$$
where $A$ stands for a uniformly distributed $\Aa$-valued random variable.

\begin{lem}\label{lem1}
The Markov chain $B_n$ is absorbed by the boundary $\partial\Aa$
in finite time, and $B_T$ is uniformly
distributed in $\partial\Aa$. Furthermore, the excursions $\Bb$ and $\Bb' \star A$ are distributed on $\Ca$ according to the same distribution $\pb : \cb\in \Ca\mapsto
\pb(\cb):= {1}/{N^{N+(|\cb|-1)}} 
$.
\end{lem}

Indeed, using the very crude upper bound
\begin{equation}\label{1est}
\PP(T>n)\leq \PP\left(\forall 0\leq p\leq n~~A_p\not\in\partial \Aa \right)=\left(1-\frac{1}{N^{N-1}}\right)^{n+1}
\end{equation}
we readily check that the chain $B_n$ is absorbed in the boundary subset $\partial\Aa$
in finite time:
$\PP(T<\infty)=1$. Furthermore, by symmetry arguments, the entrance point $B_T$ is uniformly
distributed in $\partial\Aa$, that is we have that  $\PP(B_T=a)=\nu_1(a)$. 

From the computation of the conditional expectation in Eq.~\ref{condexp}, it is easily
checked that $\pb$ is the distribution of the excursion $\Bb$. Notice that it
is a well defined probability measure on the set
of excursions $\Ca$ since
$$
\pb(\Ca_n)=\PP(T=n)\Rightarrow \pb(\Ca)=\PP(T<\infty)=1
$$
By construction, we have 
$$
\Bb'=(B'_0,\ldots,B'_{T},B'_{T+1})=(Id,\Bb)=(Id,B_0,\ldots,B_T),
$$
with $T=\{inf\ n\geq 0\ B_n\in\partial \Aa\}$.
This implies that
$$
\Bb' \star A:=(B'_0A,B'_1A,\ldots,B'_{S}A)=(A,B_0A,B_1A,\ldots,B_{S-1}A)
$$
with the stopping time
$$
S=\inf{\left\{n\geq 0~:~B'_{n} A\in\partial \Aa
\right\}}=\inf{\left\{n\geq 0~:~B_{n-1} A\in\partial \Aa
\right\}}
$$
and the convention $B_{-1}=Id$. 
The lemma follows since, by the very definition of $B_n$, the Markov chain $B_{n-1}A,\ n\geq 0$ has the same distribution as the Markov chain $B_n$. In more concrete terms, and for further use, we notice that,
for any test function $f$ we have:
\begin{eqnarray}
\EE(f(\Bb' \star A))&=&\sum_{n\geq 0}\EE\left(
f(B_{-1}A,B_0A,B_1A,\ldots,B_{n-1}A)
~1_{S=n}\right)\nonumber\\
&=&\sum_{n\geq 0}
\sum_{(c_0,\ldots,c_n)\in\Cb_{n}}f(c_0,\ldots,c_n) 
\nu(c_{0})K(c_{0},c_{1})\ldots K(c_{n-1},c_{n})
\nonumber\\
&=&\sum_{\cb\in\Cb}f(\cb)~\pb(\cb)
=\EE(f(\Bb))\label{refexc}
\end{eqnarray}
(where $\nu$ stand for the uniform probability measure on $\mathcal{A}$).
\cqfd

\begin{defi}
We associate with a probability measure $\eta\in\Pa(E)$, and  a weakly decreasing sequence of mappings  $\bb=(b_0,\ldots,b_n)\in\Aa^{n+1}$,
 a 
probability measure  $\eta_{(b_0,\ldots,b_n)}\in\Pa(E^N)$ defined for any $F\in\Ba_b(E^N)$
by
$$
\eta_{(b_0,\ldots,b_n)}(F):=\EE_{\eta}\left(F\left(\Xa_{n}^{b_n,\ldots,b_0}\right)\right)
$$
\end{defi}
By the definition
of the neutral genetic model,
for any weakly decreasing sequence of mappings  $\bb$ in
$\Aa$, and any mapping $a\in\Aa$ we have
\begin{equation}\label{key}
\eta_{\bb}\Ma=\eta_{(Id,\bb)}
\quad\mbox{and}\quad
\eta_{\bb}D_a=\eta_{\bb a}
\end{equation}
In addition, for any excursion $\cb\in\Ca$, any path $\ab=(a_1,\ldots,a_m)\in\partial \Aa^{m}$, and any $a\in\Aa$
we have
$$
\eta_{(\cb,\ab)}=(\eta M^{m})_{\cb}\qquad
\mu_{(\cb,\ab)}=\mu_{\cb}\quad\mbox{and therefore}\quad
\mu_{\cb a}=\mu_{\cb\star a}
$$
Since $(A_0,\ldots,A_n)$ has the same distribution as the reversed
sequence $(A_n,\ldots,A_0)$, we  also find that
\begin{eqnarray*}
\widehat{\Gamma}_{\eta,n}(F)&=&\EE_{\eta}\left(F\left(\Xa_n^{A_{0,n},A_{1,n},\ldots,A_{n-1,n},A_n}\right)\right)\\
&=&\EE_{\eta}\left(F\left(\Xa_n^{A_{n}\ldots A_0,A_{n-1}\ldots A_0,\ldots,A_1A_0,A_0}\right)\right)\\
&=&\EE_{\eta}\left(F\left(\Xa_n^{B_{n},B_{n-1},\ldots,B_1,B_0}\right)\right)=\EE\left(\eta_{(B_0,\ldots,B_n)}(F)\right)
\end{eqnarray*}

\begin{prop}\label{conv}
If the Markov transition $M$ has an invariant probability measure
$\mu$ on $E$, then the probability measure 
$$
\widehat{\Gamma}_{\mu}(F):=\sum_{\cb\in\Cb}
\pb(\cb)~\mu_{\cb}(F)=\EE\left(\mu_{(B_0,\ldots,B_T)}(F)\right)
$$ is an invariant measure
of the neutral genetic model $(\widehat{\xi}_n)_{n\geq 0}$.
Under the regularity condition {\rm\bf (H)}, the neutral genetic model has a unique invariant measure $\widehat{\Gamma}_{\mu}$,  and we have the estimate
\begin{equation}\label{inter}
\forall n\geq 0\qquad\left\|\widehat{\Gamma}_{\eta,n}-\widehat{\Gamma}_{\mu}\right\|_{\rm tv}\leq 
\delta~\EE(e^{-\lambda (n-T)}~1_{T\leq n})+2~\PP(T>n)
\end{equation}
with the pair of parameters $(\delta,\lambda)$ introduced in (\ref{dob}).
\end{prop}  

Indeed, if we take the excursion
$\Bb=(B_0,B_1,\ldots,B_T)$, then by (\ref{refexc}) we find that
\begin{eqnarray*}
\EE\left(\mu_{\Bb}M^{\otimes N}D_A(F)\right)&=&
\EE\left(\mu_{(Id,\Bb)}D_A(F)\right)\\
&=&\EE\left(\mu_{(Id,\Bb) A}(F)\right)=
\EE\left(\mu_{(Id,\Bb) \star A}(F)\right)=
\EE\left(\mu_{\Bb}(F)\right)
\end{eqnarray*}
Since $\EE\left(\mu_{B}(F)\right)=\widehat{\Gamma}_{\mu}(F)$, the end of the proof of the first assertion of the theorem is completed.
To prepare the proof of the second assertion, firstly we notice that
$$
\eta_{(\cb,a)}(F)=\EE_{\eta}\left(F\left(\Xa_{n+1}^{(a,c_n,\ldots,c_0)}\right)\right)=\EE_{\eta M}\left(F\left(\Xa_{n}^{(c_n,\ldots,c_0)}\right)\right)=(\eta M)_{\cb}(F)
$$
for any excursion $\cb=(c_0,\ldots,c_n)\in\Ca$, with lenght $l(\cb)=n$, 
and any $a\in\partial \Aa$. More generally, for any path $\ab=(a_1,\ldots,a_m)\in\partial \Aa^{m}$,
we have
$$
\eta_{(\cb,\ab)}=(\eta M^{m})_{\cb}
$$
This clearly implies that
\begin{eqnarray*}
\widehat{\Gamma}_{\eta,n}(F)&=&
\EE\left(\eta_{(B_0,\ldots,B_n)}(F)~1_{T\leq n}\right)+\EE\left(\eta_{(B_0,\ldots,B_n)}(F)~1_{T> n}\right)\\
&=&\EE\left((\eta M^{n-T})_{(B_0,\ldots,B_T)}(F)~1_{T\leq n}\right)+
\EE\left(\eta_{(B_0,\ldots,B_n)}(F)~1_{T> n}\right)
\end{eqnarray*}
To take the final step, we consider the decomposition
\begin{eqnarray*}
\widehat{\Gamma}_{\eta,n}(F)-\widehat{\Gamma}_{\mu}(F)&=&
\EE\left(\left[(\eta M^{n-T})_{(B_0,\ldots,B_T)}-\mu_{(B_0,\ldots,B_T)}\right](F)~1_{T\leq n}\right)\\
&&\qquad\qquad\qquad+
\EE\left(\eta_{(B_0,\ldots,B_n)}(F)-\mu_{(B_0,\ldots,B_n)}(F))~1_{T> n}\right)
\end{eqnarray*}
For any excursion $\cb=(c_0,\ldots,c_n)\in\Ca$, with lenght $l(\cb)=n$, 
we notice that
$$
[\eta_{\cb}-\mu_{\cb}](F)=\int_E~[\eta-\mu](dx)\times
\EE_{\delta_x}\left(F\left(\Xa_{n}^{(c_n,\ldots,c_0)}\right)\right)
$$
Since $\mu =\mu M^{n-T}$, applying the majoration \ref{dob}, we get
$$
\left|\widehat{\Gamma}_{\eta,n}(F)-\widehat{\Gamma}_{\mu}(F)\right|\leq 
\delta~\EE(e^{-\lambda (n-T)}~1_{T\leq n})+2\PP(T>n)
$$
for any $\|F\|\leq 1$.
This ends the proof of the proposition.
\cqfd

\section{Absorption times behavior}

We study here the renormalized 
time $T/N$ for large integers $N$, where $T$ is defined in (\ref{defT}), when
the mapping-valued Markov chain $(B_n)_{n\in\NN}$ has the identity as initial state.
Let us recall that $T$ can be interpreted as the time a neutral genetic model
with $N$ particles has to look backward to encounter its first common ancestor.
The main result is the convergence in law of $T/N$, which shows that $T$
is of order $N$, but we will also be interested in more quantitative bounds in this direction.\par\medskip
To begin with, we notice that 
$T$ only depends on $(\vert B_n\vert)_{n\in\NN}$ and that this $[N]$-valued stochastic chain is Markovian,
with transitions described by 
\begin{equation*}
\fo n\in\NN,\,\fo 1\leq p\leq q\leq N,\qquad
\PP(\vert B_{n+1}\vert=p~|~\vert B_n\vert=q)=S(q,p)\frac{(N)_p}{N^q}
\end{equation*}
(where $S(q,p)$ is the Stirling number of the second kind
giving   the number of ways of partitioning the set $[q]$ into
 $p$ non empty blocks and where  $(N)_p:={N!}/{(N-p)!}$) and starting from  $N$ if $B_0=\Id$.\par
This observation leads us to define for $N\in \NN^*$, a triangular transition matrix $M^{(N)}$ by
\bq
\fo 1\leq p, q\leq N,\qquad M_{q,p}^{(N)}&\df& S(q,p)\frac{(N)_p}{N^q}\eq
and to consider for any $1\leq i\leq N$, a Markov chain $R^{(N,i)}\df (R^{(N,i)}_n)_{n\in\NN}$
starting from $i$ and whose transitions are governed by $M^{(N)}$. Such a Markov chain
is (a.s.) non-increasing and 1 is an absorption state. We denote
\bq
S^{(N,i)}&\df&\inf\{n\in\NN\st R^{(N,i)}_n=1\}\eq
so that $T$ has the same law as $S^{(N,N)}$. The goal of this section is to prove the
\begin{theo}\label{th1}
Let $(i_N)_{N\in\NN^*}$ be a sequence of integers satisfying $1\leq i_N\leq N$ for any $N\in\NN^*$
and diverging to infinity. Then the following convergence in law takes place for large $N$
\bq
\frac{S^{(N,i_N)}}{N}&\overset{\La}{\longrightarrow}& \sum_{l\in\NN}\frac{2}{(l+1)(l+2)}\Ea_l\eq
where $(\Ea_l)_{l\in\NN}$ is an independent family of exponential variables of parameter 1.
Furthermore, we have for any fixed $0\leq \alpha<1$,
\bq
\sup_{N\in\NN^*}\EE[\exp(\alpha S^{(N,i_N)}/N)]&\leq  & \exp(\alpha)\Pi(\alpha)\eq
where 
\bq
\fo 0\leq \alpha<1,\qquad \Pi(\alpha)&\df& \prod_{l\in\NN}\frac1{1-\frac{2\alpha}{(l+1)(l+2)}}\eq
is the Laplace transform of the above limit law.
Thus for any continuous function $f\st\RR_+\ri\RR_+$ verifying
$\lim_{\alpha\ri1_-}\limsup_{s\ri+\iy}\exp(-\alpha s)f(s)=0$, we are insured of
\bq
\lim_{N\ri\iy}\EE[f(S^{(N,i_N)}/N)]&=&\EE\lt[f\lt( \sum_{l\in\NN}\frac{2}{(l+1)(l+2)}\Ea_l\rt)\rt]\eq
\end{theo}
In particular, for any given $0\leq \alpha<1$, via a Markov inequality, we deduce
the exponential upper bound
\bq
\fo N\in\NN^*,\,\fo 1\leq i\leq N,\,\fo n\in \NN,\qquad
\PP[S^{(N,i)}\geq n]&\leq & \exp(\alpha)\Pi(\alpha)\exp(-\alpha n/N)\\
&\leq &\frac{K}{1-\alpha}\exp(-\alpha n/N)
\eq
with $K=e\left(\prod_{l\in\NN^*}\left(1-\frac{2}{(l+1)(l+2)}\right)\right)^{-1}$. Optimizing 
this inequality with respect to $0\leq \alpha<1$, we obtain that for any positive integers $n$ and $N$,
 
\bqn{NeN}
\fo 1\leq i\leq N,\qquad
\PP[S^{(N,i)}\geq n]&\leq & K \left(\frac{n}{N}\vee 1\right)\exp(-( n/N-1)_+) \label{majtemps}\eqn
(of course this bound does not give any relevant information for $n\leq N$, when 
the l.h.s.\ probability is not small).
As it was explained in the second section, such an inequality is useful to 
estimate convergence to equilibrium for neutral genetic models
and the results presented in the introduction follow from it.
\\
Indeed, in view of proposition \ref{conv}, we still need another estimate, but it is an immediate
consequence of the above bound:
\begin{cor} \label{majexp}
There exists a constant $K'\geq 0$ such that for any $\lambda> 0$,
any $N\in\NN^*$, any $1\leq i\leq N$ and any $n\geq N+1/\lambda$, we have
\bq
\EE\lt[\exp(-\lambda(n-S^{(N,i)}))\un_{\{S^{(N,i)}\leq n\}}\rt]&\leq & K'\frac{n}{N}\exp\lt(-\lt(1+\frac1{\lambda N}\rt)^{-1}\frac{n}{N}\rt)\eq
\end{cor}

Indeed, let $m\in[n]$ be given, we have
\begin{eqnarray*}
\lefteqn{\EE\left[\exp({-\lambda (n-S^{(N,i)})})\un_{\{S^{(N,i)}\leq n\}}\right]}\\&=&\EE\left[\exp({-\lambda(n-S^{(N,i)})})\un_{\{S^{(N,i)}\leq m\}}\right]+\EE\left[\exp({-\lambda(n-S^{(N,i)})})\un_{\{m\leq S^{(N,i)}\leq n\}}\right]\\
&\leq &\exp({-\lambda(n-m)})+\PP(S^{(N,i)}>m)\\
&\leq &\exp({-\lambda(n-m)})+K \left(\frac{m}{N}\vee 1\right)\exp(-( m/N-1)_+)\\
&\leq &(1+K)\frac{n}{N}\max(\exp({-\lambda(n-m)}),\exp(-( m/N-1)_+))
\end{eqnarray*}
Where we have used that $m\leq n$ and that $n\geq N$. Optimizing the last term, we are led to consider
\bq
m&=&\lt\lfloor\frac{n}{1+\frac1{\lambda N}}\rt\rfloor\eq
Note that this integer number is larger than $N$ for $n\geq N+1/\lambda$
and the corollary's inequality follows easily from
the obvious bounds
\bq
\frac{n}{1+\frac1{\lambda N}}-1\ \leq \ m\ \leq \ \frac{n}{1+\frac1{\lambda N}}\eq
\wwtbp
\begin{rem}
The limit distribution appearing in Theorem~\ref{th1} is the same
as the law of the coalescence time for the Kingman process (see for instance \cite{kingman2}).
This could have been expected, since it is known that, suitably ``rearranged'', 
the mapping-valued Markov process $(B_{\lfloor Nt\rfloor})_{t\in\RR_+}$ converges to the Kingman coalescent process for large $N$.
Nevertheless, at our best knowledge, this convergence takes place in a weak sense
which does not permit to deduce the results presented here. In fact, the latter could serve to strengthen the previous
convergence.

\end{rem}
Before proving Theorem \ref{th1}, we will investigate the simpler
problem where only negative jumps of unit length are permitted. More precisely,
let $\wi M^{(N)}$ be the transition matrix defined by
\bq
\fo 1\leq q,p\leq N,\qquad \wi M^{(N)}_{q,p}&=&\lt\{\begin{array}{ll}
\di\frac{(N)_q}{N^q} &\hbox{, if $p=q$}\\
\di 1-\frac{(N)_q}{N^q}&\hbox{, if $p=q-1$}\\
\di 0&\hbox{, otherwise}\end{array}\rt.
\eq
Every corresponding notion will be overlined by a tilde. 
Then we have a result which is similar 
to Theorem \ref{th1}, with nevertheless some slight differences:
\begin{prop}\label{pr1}
Let $(i_N)_{N\in\NN^*}$ be a sequence of integers satisfying $1\leq i_N\leq N$ for any $N\in\NN^*$,
diverging to infinity
and such that $a\df \lim_{N\ri\iy}i_N/N$ exists in $[0,1]$. Then the following convergence in law takes place for large $N$
\bq
\frac{\wi S^{(N,i_N)}}{N}&\overset{\La}{\longrightarrow}& a+\sum_{l\in\NN}\frac{2}{(l+1)(l+2)}\Ea_l\eq
where $(\Ea_l)_{l\in\NN}$ is an independent family of exponential variables of parameter 1.
Furthermore, we have for any fixed $0\leq \alpha<1$,
\bq
\sup_{N\in\NN^*}\EE\lt[\exp\lt(\alpha \frac{\wi S^{(N,i_N)}-i_N+1}{N}\rt)\rt]&=  & \Pi(\alpha)\eq
\end{prop}
The case of a fixed initial condition, i.e.\ when there exists $i\in\NN\setminus\{0,1\}$ such that for any
$N\geq i$, $i_N=i$, is also instructive, despite the fact it is not included in the previous result:

\begin{lem}\label{l1}
For given $i\in\NN\setminus\{0,1\}$,  the following convergence in law takes place for large $N$,
\bq
\frac{\wi S^{(N,i)}}{N}&\overset{\La}{\longrightarrow}& \sum_{0\leq l\leq i-2}\frac{2}{(l+1)(l+2)}\Ea_l\eq
where the $\Ea_l$, for $0\leq l\leq i-2$, are independent exponential variables of parameter 1.
Furthermore, we have for any fixed $0\leq \alpha<1$,
\bq
\sup_{N\in\NN^*, N\geq i}\EE\lt[\exp\lt(\alpha \frac{\wi S^{(N,i)}-i+1}{N}\rt)\rt]&=  & 
\prod_{0\leq l\leq i-2}\frac1{1-\frac{2\alpha}{(l+1)(l+2)}}
\eq
\end{lem}

Let $2\leq i\leq N$ be given. By a backward iteration and with the convention that $\wi T^{(N)}_{i}=0$, we define for $1\leq j\leq i-1$,
\bq
\wi T^{(N)}_{j}&\df& \inf\lt\{ n\in\NN\st \wi R^{(N,i)}_{\wi T^{(N)}_{j+1}+n}=j\rt\}\eq
In words, $\wi T^{(N)}_{j}$ is the time necessary for  the Markov chain $(\wi R^{(N,i)}_n)_{n\in\NN}$
to jump from $j+1$ to $j$. For $1\leq j\leq i-1$, let us also denote $\wit T^{(N)}_{j}\df
\wi T^{(N)}_{j}-1$ and $p_{N,j}\df \wi M^{(N)}_{j+1,j}= (N)_{j+1}/N^{j+1}$.
It is clear that  $\wit T^{(N)}_{j}$ 
is distributed as a geometric law of parameter $p_{N,j}$, namely,
\bq
\fo m\in\NN,\qquad \PP[\wit T^{(N)}_{j}=m]&=&(1-p_{N,j})p_{N,j}^m\eq
Furthermore, the variables $\wi T^{(N)}_{j}$, for $1\leq j\leq i-1$, are independent and we can write
\bqn{decomp}
\nonumber\wi S^{(N,i)}&=&\sum_{1\leq j\leq i-1}\wi T^{(N)}_{j}\\
&=&i-1+\sum_{1\leq j\leq i-1}\wit T^{(N)}_{j}\eqn
This leads us to study the individual behavior of the summands:
\begin{lem}
\label{l2}
With the above notation and for a fixed $j\in\NN^*$, we are insured  of the following convergence in law as $N$ goes to infinity,
\bq
\frac{\wit T_j^{(N)}}{N}&\overset{\La}{\longrightarrow}& \frac{2}{j(j+1)}\Ea\eq
where  $\Ea$ is an exponential variable of parameter 1.
Furthermore, we have for any fixed $0\leq \alpha<j(j+1)/2$,
\bq
\sup_{N\in\NN^*, N\geq j+1}\EE\lt[\exp\lt(\alpha \frac{\wit T_j^{(N)} }{N}\rt)\rt]&=  & 
\lim_{N\ri\iy}\EE\lt[\exp\lt(\alpha \frac{\wit T_j^{(N)} }{N}\rt)\rt]\\&=&
\frac1{1-\frac{2\alpha}{j(j+1)}}
\eq
\end{lem}

The simplest way to prove the Lemma seems to resort to Laplace's transform.  So we compute that for any $\alpha\in\RR$,
\bq
\EE[\exp(\alpha\wit T_j^{(N)}/N)]&=&\lt\{
\begin{array}{ll}
\di\frac{1-p_{N,j}}{1-\exp(\alpha/N)p_{N,j}}&\hbox{, if $\exp(\alpha/N)p_{N,j}<1$}\\
\di+\iy&\hbox{, otherwise}\end{array}\rt.\eq
The condition $\exp(\alpha/N)p_{N,j}<1$ is equivalent to 
\bq
\alpha&<&-N\sum_{1\leq k\leq j}\ln\lt(1-\frac{k}{N}\rt)\eq
and taking into account the convexity inequality $-\ln(1-x)\geq x$, valid for any $x<1$, 
we get that it is in particular fulfilled for $0\leq \alpha<j(j+1)/2$.\\
Furthermore, using an asymptotic expansion of $p_{N,j}$ in $1/N$, we show without difficulty 
that uniformly for $\alpha$ on any compact set of $(-\iy,j(j+1)/2)$ (in particular 
in some neighborhoods of 0),
\bq
\lim_{N\ri\iy}\frac{1-p_{N,j}}{1-\exp(\alpha/N)p_{N,j}}&=&\frac{1}{1-\frac{2\alpha}{j(j+1)}}\eq
 Since for $\alpha<j(j+1)/2$, the above r.h.s.\ coincides with the Laplace's transform
of  $\frac{2}{j(j+1)}\Ea$, a well-known result (see for instance 
Theorem 0.5 of \cite{toulouse}, indeed, as we are working 
with nonnegative random variables, only a right neighborhood of 0 is required) enables us to conclude to the announced convergence in law.
\\
Concerning the upper bound, we begin by remarking that by a previously mentioned convexity inequality,
we have
\bq
\fo N\geq j,\qquad p_{N,j}&\leq &\exp\lt(-\frac{j(j+1)}{2N}\rt)\eq
Next, we notice that for $\alpha\geq 0$, the mapping
\bq
[0,\exp(-\alpha/N))\ni t&\mapsto& \frac{1-t}{1-\exp(\alpha/N)t}\eq
is increasing, so that
\bq
\frac{1-p_{N,j}}{1-\exp(\alpha/N)p_{N,j}}&\leq &\frac{1-\exp\lt(-\frac{j(j+1)}{2N}\rt)}{1-\exp(\alpha/N)\exp\lt(-\frac{j(j+1)}{2N}\rt)}\eq
Finally, we consider for fixed $x>0$, the function
\bq
\varphi\st [0,x]\ni y&\mapsto& x(1-\exp(y-x))+(y-x)(1-\exp(-x))\eq
A variation study shows that it is increasing up to some point $y_x$ belonging to $(0,x)$
and decreasing after this point. Since $\varphi(0)=\varphi(x)=0$, we 
get that $\varphi$ is nonnegative on $[0,x]$.
Thus we obtain that
\bq
\fo 0\leq y<x,\qquad \frac{1-\exp(-x)}{1-\exp(y-x)}&\leq &\frac{x}{x-y}\eq
Applying this inequality with $x=j(j+1)/(2N)$ and $y=\alpha/N$, it appears that for
$0\leq \alpha<j(j+1)/2$ and $N\geq j+1$,
\bq
\EE\lt[\exp\lt(\alpha \frac{\wit T_j^{(N)} }{N}\rt)\rt]&\leq &
\frac1{1-\frac{2\alpha}{j(j+1)}}\eq
Since the l.h.s.\ was already seen to converge to the r.h.s.\ for large $N$,
we can conclude that the desired equalities hold.\wwtbp

\noindent The proof of Lemma \ref{l1} follows at once
from the  decomposition (\ref{decomp})
and Lemma \ref{l2}.
\wwtbp

The proof of Proposition \ref{pr1} is based on arguments that are close to the preceding ones.
More precisely, similarly to (\ref{decomp}),
we can write
\bq\wi S^{(N,i_N)}-i_N+1&=&\sum_{1\leq j\leq i_N-1}\wit T^{(N)}_{j}
\eq
Thus we get for any $\alpha<1$,
\bq
\EE\lt[\exp\lt(\alpha \frac{\wi S^{(N,i_N)}-i_N+1}{N}\rt)\rt]&=&
\prod_{1\leq j\leq i_N-1}\frac{1-p_{N,j}}{1-\exp(\alpha/N)p_{N,j}}\\
&\leq &\prod_{1\leq j\leq i_N-1}\frac1{1-\frac{2\alpha}{j(j+1)}}\\
&\leq &\Pi(\alpha)\eq
Let $i\in\NN^*$ be fixed (at first). By assumption, for $N$ large enough,
we will have $i_N\geq i$ and quite obviously, $\wi S^{(N,i_N)}-i_N+1$ will stochastically dominate
$\wi S^{(N,i)}-i+1$, which implies, for $0\leq \alpha<1$,
\bq
\EE\lt[\exp\lt(\alpha \frac{\wi S^{(N,i_N)}-i_N+1}{N}\rt)\rt]&\geq &
\EE\lt[\exp\lt(\alpha \frac{\wi S^{(N,i)}-i+1}{N}\rt)\rt]\eq
and thus
\bq
\liminf_{N\ri\iy}\EE\lt[\exp\lt(\alpha \frac{\wi S^{(N,i_N)}-i_N+1}{N}\rt)\rt]
&\geq &\lim_{N\ri\iy}\EE\lt[\exp\lt(\alpha \frac{\wi S^{(N,i)}-i+1}{N}\rt)\rt]\\
&=&\prod_{1\leq j\leq i-1}\frac1{1-\frac{2\alpha}{j(j+1)}}\eq
Letting $i$ go to infinity, it appears that
\bq
\lim_{N\ri\iy}\EE\lt[\exp\lt(\alpha \frac{\wi S^{(N,i_N)}-i_N+1}{N}\rt)\rt]
&=&\Pi(\alpha)\eq
and the last part of Proposition \ref{pr1} follows.
We also obtain that for $0\leq \alpha<1$,
\bq
\lim_{N\ri\iy}\EE\lt[\exp\lt(\alpha \frac{\wi S^{(N,i_N)}}{N}\rt)\rt]
&=&\exp(\alpha a)\Pi(\alpha)\eq
and to conclude to the announced convergence in law,
it remains to check that above convergence is uniform in some compact
right neighborhood of 0.
But this is a consequence of Dini's theorem, since
the r.h.s.\ is continuous in $\alpha\in[0,1)$
and  for all fixed $N\in\NN^*$, the mapping $[0,1)\ni\alpha\mapsto
\EE[\exp(\alpha {\wi S^{(N,i_N)}}/{N})]$
is increasing.
\par\me
We now proceed to the proof of Theorem \ref{th1}.
Its second part is the simplest one, since $S^{(N,i_N)}$ is stochastically
dominated by $\wi S^{(N,i_N)}$:
\bq
\fo n\in \NN,\qquad
\PP[S^{(N,i_N)}\geq n]&\leq &\PP[\wi S^{(N,i_N)}\geq n]\eq
This fact is based on two observations: on one hand, assuming that the Markov chain $R^{(N,i_N)}$
is in state $j\in[i_N]$ at some time $n$, it will wait the same time
to jump out of it as $\wi R^{(N,i_N)}$, but $R^{(N,i_N)}$
will visit less states in $[i_N]$ than $\wi R^{(N,i_N)}$.
Taking into account these two facts, the above inequalities follow immediately. Details are left to the reader: one can e.g. construct a coupling between
$R^{(N,i_N)}$ and $\wi R^{(N,i_N)}$  such that
$\PP[\wi S^{(N,i_N)}\geq S^{(N,i_N)}]=1$ (for a general reference on the subject, see for instance \cite{lindvall}).
\\
In particular we get that for any $\alpha\geq 0$,
\bq
\sup_{N\in\NN^*}\EE[\exp(\alpha S^{(N,i_N)}/N)]&\leq  &\sup_{N\in\NN^*}\EE[\exp(\alpha \wi S^{(N,i_N)}/N)]\eq
so the wanted upper bound follows from that of Proposition \ref{pr1}.
\par
We will need to work more to obtain the convergence in law. Heuristically the proof is based on the fact
\begin{itemize}
\item that the Markov chain $R^{(N,i_N)}$ will rapidly reach some point negligible with respect to $N$
(but however going to infinity with $N$) 
\item and that from this point to 1, $R^{(N,i_N)}$
and
$\wi R^{(N,i_N)}$ are quite similar, which will enable us to make use of Proposition \ref{pr1}.
\end{itemize}
\par
\sm

We begin by showing  the second assertion, namely that for not too large $i_N$,
$R^{(N,i_N)}$ and $\wi R^{(N,i_N)}$ are almost the same.
The following lemma will enable us to quantify the ``not too large''.

\begin{nota}
 In the following, we will sometimes drop the superscript $(N,i)$ in $R^{(N,i)}$ and $\tau^{(N,i)}$ (to be defined below) when no confusion is possible, in order to make the proofs more readable.
\end{nota}

\begin{lem}\label{sautdedeux}
There exists a constant $\chi_1>0$ such that for any $2\leq q\leq N/2$,
\bq
\fo n\in\NN,\fo i \in \NN^*,\qquad \PP[R_{n+1}^{(N,i)}\leq R_n^{(N,i)}-2\vert R_n^{(N,i)}=q]&\leq & \chi_1\frac{q^4}{N^2}\eq
\end{lem}

Indeed, by definition, we have for any $2\leq q\leq N$ and $n\in\NN$,
\bq
\PP[R_{n+1}\leq R_n-2\vert R_n=q]&=&1-\PP[R_{n+1}= R_n\vert R_n=q]-\PP[R_{n+1}= R_n-1\vert R_n=q]\\
&=&1-\frac1{N^q}(S(q,q)(N)_q+S(q,q-1)(N)_{q-1})
\\
&=&1-\frac{(N)_{q-1}}{N^q}\lt(N-q+1+\frac{q(q-1)}{2}\rt)\\
&=&1-\exp\lt(\sum_{1\leq l\leq q-2}\ln(1-l/N)\rt)\lt(1+\frac{(q-1)(q-2)}{2N}\rt)\eq
But we note that there exists a constant $\chi>0$ such that
\bq
\fo 0\leq x\leq 1/2,\qquad \ln(1-x)&\geq &-x-\chi x^2\eq
thus we get
\bq
\sum_{1\leq l\leq q-2}\ln(1-l/N)&\geq & -\sum_{1\leq l\leq q-2}\left(\frac{l}{N}+\chi\frac{l^2}{N^2}\right)\\
&=&-\frac{(q-2)(q-1)}{2N}-\frac{\chi}{6N^2}(q-2)(q-1)(2q-3)\\
&\geq &-\frac{(q-2)(q-1)}{2N}-\frac{\chi}{3N^2}q^3\eq
Taking into account the convexity inequality $\exp(x)\geq 1+x$, valid for all $x\in\RR$,
the wanted probability is then bounded above by
\bq
1-\lt(1-\frac{(q-2)(q-1)}{2N}-\frac{\chi}{3N^2}q^3\rt)\lt(1+\frac{(q-1)(q-2)}{2N}\rt)\eq
expression which is dominated by $\chi_1\frac{q^4}{N^2}$ for an appropriate choice of the constant
$\chi_1$.\wwtbp
% WRONG !!!
% \begin{rem} \label{couplage}
% It follows the preceding Lemma that  the probability that $R^{(N,i_N)}$ admits at least a jump downward  strictly larger (in absolute value) than one, 
% is bounded above by $\chi_3 i^5_N/N^2$ for a well-chosen constant $\chi_3$ independent of $2\leq i_N\leq N/2$.
% \end{rem}
Let us remark that the previous lemma enables us to get useful estimates on natural couplings between 
 $R^{(N,i_N)}$ and $\wi R^{(N,i_N)}$, when $i_N$ is quite small in comparison with $N$.
More precisely, we begin by constructing
$R^{(N,i_N)}$ and we define
\bq
\wi\tau^{(N,i_N)}\df\inf\{n\geq 0\st R_{n+1}^{(N,i_N)}\leq R_{n}^{(N,i_N)}-2\}\eq
Then we take $\wi R_{n}^{(N,i_N)}=R_{n}^{(N,i_N)}$ for $n\leq \wi\tau^{(N,i_N)}$
and $\wi R_{\wi\tau^{(N,i_N)}+1}^{(N,i_N)}=R_{\wi\tau^{(N,i_N)}}^{(N,i_N)}-1$.
Next,  we use the transition matrix $\wi M^{(N)}$
to construct $(\wi R_{n}^{(N,i_N)})_{n>\wi\tau^{(N,i_N)}}$ in a traditional way
(and independently from the above constructions).
It is easy to check that this gives a coupling between $R^{(N,i_N)}$ and $\wi R^{(N,i_N)}$
and we have for any $A_N>0$, at least if $i_N\leq N/2$,
\bq
\lefteqn{\PP[\ex n\in\NN\st R^{(N,i_N)}_n\not=\wi R^{(N,i_N)}_n]}\\&= &
\PP[\wi\tau^{(N,i_N)}<S(N,i_N)]\\
&\leq&\PP[\wi\tau^{(N,i_N)}<S(N,i_N),\,S(N,i_N)\leq A_NN]+ \PP[S(N,i_N)>A_NN]\\
&\leq &A_NN\max_{2\leq q\leq i_N}\chi_1\frac{q^4}{N^2}+\PP[\wi S(N,i_N)>A_NN]\\
&=&\chi_1A_N\frac{i_N^4}{N}+\PP[\wi S(N,i_N)>A_NN]\eq
Thus it appears that if the sequence $(i_N)_{N\in\NN^*}$ of initial states diverging to infinity satisfies furthermore
\bqn{iN5}
\lim_{N\ri\iy}\frac{i^4_N}{N}&=&0\eqn
then we can choose a sequence $(A_N)_{N\in\NN^*}$
diverging to infinity and verifying 
\bq 
\lim_{N\ri\iy} A_N\frac{i_N^4}{N}&=&0\eq 
so that
\bq 
\lim_{N\ri\iy} \PP[ R^{(N,i_N)}\not=\wi R^{(N,i_N)}]&=&0\eq 
In particular
\bq
\lim_{N\ri\iy}\PP[S^{(N,i_N)}\not=\wi S^{(N,i_N)}]&=&0\eq
and we get 
from proposition \ref{pr1}
(applied with $a=0$, because (\ref{iN5}) implies that $\lim_{N\ri\iy}i_N/N=0$)
that
\bq
\frac{S^{(N,i_N)}}{N}&\overset{\La}{\longrightarrow}& \sum_{l\in\NN}\frac{2}{(l+1)(l+2)}\Ea_l\eq
\par
To treat the general case,
we now turn our attention toward the first assertion 
before Lemma \ref{sautdedeux} (that the Markov chain $R^{(N,i_N)}$ will rapidly reach some point negligible with respect to $N$). So for $1\leq j\leq i\leq N$, we consider the reaching time
\bq
\tau_j^{(N,i)}&\df&\inf\{n\in\NN\st R^{(N,i)}_n<j\}\eq
It will be very useful to know
 that when $R^{(N,i)}$ is jumping through $j$,
it is not going too far away from $j$.
The next lemma gives an estimate in this direction:
\begin{lem}\label{saut}
There exists a constant $\chi_2>0$, such that for all $4\leq j\leq \sqrt{N}$
and all $j\leq i\leq N$, we have
\bq
\PP\lt[R^{(N,i)}_{\tau_j^{(N,i)}}\leq j/2\rt]&\leq &\chi_2\frac{j^4}{N^{2}}\eq
\end{lem}
Let us prove the Lemma. For simplicity, we write $k\df\lfloor j/2\rfloor$, and for $1\leq l\leq k$, $j\leq m\leq N$ and
$n\in\NN$,
we consider the quantity
\bq
\frac{\PP[R_{n+1}=l\vert R_n=m]}{\PP[R_{n+1}=l+k\vert R_n=m]}
&=&\frac{S(m,l)}{S(m,l+k)}\frac{(N)_l}{(N)_{k+l}}\eq
We remark that $S(m,l)\leq S(m,l+k)S(l+k,l)$, because if we have a partitioning of $[m]$
into $l+k$ blocks (that we can order through their respective smaller elements) and a partitioning of $[l+k]$ into $l$ blocks,
we can naturally construct a partitioning of $[m]$ into
$l$ blocks by composing them and this mapping is clearly onto.
Thus, since $S(l+k,l+k)=1$, we have
\bq
\frac{\PP[R_{n+1}=l\vert R_n=m]}{\PP[R_{n+1}=l+k\vert R_n=m]}
&\leq &\frac{\PP[R_{n+1}=l\vert R_n=l+k]}{\PP[R_{n+1}=l+k\vert R_n=l+k]}\eq
But
the last denominator is 
\bq
\frac{(N)_{l+k}}{N_{l+k}}&=&\exp\lt(\sum_{1\leq q\leq l+k-1}\ln\lt(1-\frac{q}{N}\rt)\rt)
\eq
expression which is bounded below by a positive constant, uniformly in $N\in\NN$ and
$l+k\leq j\leq \sqrt{N}$ (as it was seen in the proof of lemma \ref{sautdedeux}). Since $k\geq 2$, we have
\bq
\PP[R_{n+1}=l\vert R_n=l+k]&\leq &\PP[R_{n+1}\leq l+k-2\vert R_n=l+k]\eq
and we have seen in lemma \ref{facile} that this quantity is bounded above
by $j^4/N^2$ up to an universal constant, for $l+k\leq j$.
Thus there exists a constant $\chi>0$ such that for any $N,k,l,m,n$ as above, 
we have (since for $4\leq j\leq \sqrt{N}$, we have $\sqrt{N}\leq N/2$),
\bq
{\PP[R_{n+1}=l\vert R_n=m]}\leq \chi\frac{j^4}{N^{2}}\PP[R_{n+1}=l+k\vert R_n=m]\eq
and summing these inequalities for $1\leq l\leq k$, we get\bq
\PP[R_{n+1}\leq k\vert R_n=m]&\leq &\chi\frac{j^4}{N^{2}}\PP[k<R_{n+1}\leq 2k\vert R_n=m]\\
&\leq &\chi\frac{j^4}{N^{2}}\PP[R_{n+1}\leq 2k\vert R_n=m]\eq
This bound can also be rewriten
\bq
\PP[R_{n+1}\leq k,  R_n=m,\tau_j=n]&\leq &\chi\frac{j^4}{N^{2}}\PP[ R_n=m,\tau_j=n]\eq
and summing over all $n\in\NN$ and
$j\leq m\leq N$, we obtain
finally the wanted inequality.\wwtbp
With this estimate at hand, we will be able
to take advantage of an approach inspired by section 1.4 of the book \cite{motwani} of Motwani and Raghavan,
to investigate the expectation of $\tau_j^{(N,i)}$.
\begin{lem}\label{MR}
There exists a constant $\chi_3>0$ such that for all $4\leq j\leq i\leq N$ with $j\leq N^{2/5}$, we have
\bq
\EE[\tau_j^{(N,i)}]&\leq &\chi_3 \frac{N}{j}\eq
\end{lem}
Thus the Markov chain $R^{(N,N)}$ goes relatively fast from $N$ to $\lfloor N^{2/5} \rfloor$.\par\me

To prove the Lemma, notice first that, by definition of the chain $R$, we find that for any $n\in\NN$ and any $1\leq q\leq N$,
\bq
\EE\left[N-R_{n+1}~|~R_{n}=q\right]&=&\frac{1}{N^q}\sum_{1\leq p\leq q}~S(q,p)~(N)_p~(N-p)
\eq
Using the fact that
$
(N)_p~(N-p)=N~(N-1)_p$,
we find that
\bq
\EE\left[N-R_{n+1}~|~R_{n}=q\right]&=&\frac{1}{N^{q-1}}\sum_{1\leq p\leq q}~S(q,p)~(N-1)_p\\
&=& \frac{(N-1)^q}{N^{q-1}}\\&=&N~\left(1-\frac{1}{N}\right)^q\eq
This implies that for any $1\leq q\leq N$,
\bq
\EE\left[R_{n}-R_{n+1}~|~R_{n}=q\right]&=&\EE\left[q-N+N-R_{n+1}~|~R_{n}=q\right]\\
&=&q-N+N~\left(1-\frac{1}{N}\right)^q\eq
We are thus led to study the function defined by
\bq
\fo s\geq 1, \qquad g(s)&\df& s-N+N\left(1-\frac{1}{N}\right)^s\eq
and will show that there exists a constant $\chi$
such that
\bqn{g}
\fo N\geq 2,\,\fo 1\leq s\leq N,\qquad g(s)&\geq &\frac{(s-1)^2}{\chi N}\eqn
Indeed, we compute that
\bq
\fo s\geq 1,\qquad g''(s)&=&N\ln^2\lt(1-\frac{1}{N}\right)\left(1-\frac{1}{N}\right)^s\eq
and we see without difficulty that there exists a universal positive constant such that the r.h.s.\ is bounded below
by $1/(\chi N)$ for $N\geq 2$ and $1\leq s\leq N$.
Furthermore, we have
\bq
g'(1)&=&1+(N-1)\ln\lt(1-\frac{1}{N}\right)\eq
and as a function of $N\in[2,+\iy[$, this quantity is decreasing,
so that it is positive since its limit in $+\iy$ is 0.\\
Now, the lower bound (\ref{g}) follows from a second order Taylor-Lagrange formula.
As a by-product of the above facts,
we deduce that $g$ is increasing on $[1,+\iy)$, since $g'(1)> 0$ and $g''>0$.
Furthermore, we note that a reverse bound is valid: 
% there exists a constant $\chi>0$ such that
\bqn{facile}
\fo N\geq 2,\,\fo s\geq 0,\qquad g(s)&\leq & \frac{s^2}{2N}\eqn
Indeed, we have for any $N\geq 2$ and $r>0$, 
\bq
g(rN)&=&N\lt(r-1+\lt(1-\frac1N\rt)^{rN}\rt)\\&\leq &
N\lt(r-1+\exp(-r)\rt)\eq
but for $r> 0$, it is easily seen that the expression between parentheses
is less than $r^2/2$, so we get (\ref{facile}) by replacing $s$ by $rN$. It also follows from
(\ref{g}) and (\ref{facile})
that there
exists a constant $\wi \chi>0$ such that $g(l+1)\leq \wi\chi g(l)$ holds for any $l\geq 2$.\\
% as an easy consequence of $g(1)=0$ and of the two following inequalities,
% \bq
% g'(1)&\leq & \frac1N\\
% \fo N\geq 2,\,\fo s\geq 1,\qquad g''(s)&\leq &\frac{\chi}{N}\eq
Next we consider the stochastic chain $M$ defined iteratively by

$$
\fo n\in\NN,\qquad M_n \df n+\sum_{1\leq l\leq R_n}\frac1{g(l)}
$$

Let us check that it is a supermartingale with respect to the filtration $(\Fa_n)_{n\in\NN}$
naturally generated by the Markov chain $R$.
It is clearly adapted, so for $n\in\NN$, we compute, via the Markov property and the fact that $g$ is increasing, that
\bq
\lefteqn{\EE[M_{n+1}-M_n\vert \Fa_n]}\\
&=&
\EE\lt[1-\sum_{R_{n+1}<l\leq R_n }\frac1{g(l)}
\Bigg\vert \Fa_n\rt]\\
&\leq &\EE\lt[1-\frac{R_n-R_{n+1}}{g(R_n)}
\Bigg\vert \Fa_n\rt]
\eq
But the  last rhs is zero by definition of $g$, since it can also be written
\bq
1-\frac{\EE\lt[{R_n-R_{n+1}}
\Bigg\vert R_n \rt]}{g(R_n)}&=&0\eq

Next we apply Doob's stopping theorem to the nonnegative surmartingale $M$
with respect to the stopping time $\tau_j$, to get 
\bq
\EE[M_{\tau_j}]&\leq &\EE[M_0]\eq
But we note that on one hand,
\bq
\EE[M_0]&=&\sum_{1\leq l\leq i}\frac{1}{g(l)}\eq
and on the other hand, by definition,
\bq
M_{\tau_j}&=&\tau_j+\sum_{1\leq l\leq R_{\tau_j}}\frac{1}{g(l)}
\eq
so it appears that
\begin{eqnarray*}
\EE(\tau_j)& \leq &\EE\left(  \sum_{ R_{\tau_j} < l \leq i }  \frac{1}{g(l)}\right)\\
&\leq& \wi\chi\,\EE\left(  \sum_{ R_{\tau_j} < l \leq i }  \frac{1}{g(l+1)}\right)\\
& \leq &  \wi\chi\,\EE \left ( \int_{R_{\tau_j} }^{i} \frac{1}{g(s+1)} ds \right)\\
& \leq & \wi\chi\,\EE \left( \chi \int_{R_{\tau_j} }^i \frac{N}{s^2} ds \right)\\
& \leq & \wi\chi\chi N \EE\left( \frac{1}{R_{\tau_j} } \right)
\end{eqnarray*}
Thus it remains to evaluate the last expectation.
To do that, we decompose it into two parts:
\bq
\EE\left( \frac{1}{R_{\tau_j} } \right)
&=&\EE\left( \frac{1}{R_{\tau_j} } \un_{R_{\tau_j}\leq j/2}\right)
+\EE\left( \frac{1}{R_{\tau_j} } \un_{R_{\tau_j}> j/2}\right)
\eq
The last term is bounded above by $2/j$, while by lemma \ref{MR}, the first
term is of order $j^4/N^2$, up to an universal constant. 
Since by our assumption we have $j^4/N^2\leq 1/j$,
the annouced estimate follows at once.
\wwtbp
Now we can continue the proof of Theorem \ref{th1}.
let us consider
\bq\fo N\in\NN^*,\qquad j_N&\df&i_N\wedge \lfloor N^{1/5}\rfloor
\eq
and to simplify notations, let us write $\tau_N\df\tau^{(N,i_N)}_{j_N}$ if $j_N<i_N$ and $\tau_N=0$ if $j_N=i_N$.
Then we can decompose $S^{(N,i_N)}$ into
$\tau_N+\wi\tau_N$,
where
\bq
\wi\tau_N&\df&\inf\{n\in\NN\st R^{(N,i_N)}_{\tau_N+n}=1\}\eq
From lemma \ref{MR}, we see that $\tau_N/N$ converges in probability to zero
for large $N$. 
Thus we are led to study the behavior of $\wi\tau_N/N$ for large $N$, and we
begin by noting that
conditionally to $R^{(N,i_N)}_{\tau_N}$, $\wi \tau_N$ has the same law as 
$S^{(N,R^{(N,i_N)}_{\tau_N})}$. 
So we have to check that $R^{(N,i_N)}_{\tau_N}$ is not too small and this is insured by lemma
\ref{saut}.
\\
To conclude,
we remark that by previous estimates, we already know that the family of the distributions
of the $S^{(N,i_N)}/N$ (or equivalently of the $\wi \tau_N/N$), for $N\in\NN^*$, is relatively compact for
the weak convergence. So we just have to verify that for any converging subsequence, the limit 
coincides with the law of $\sum_{l\in\NN}\frac{2}{(l+1)(l+2)}\Ea_l$.
For this purpose, we can furthermore assume, up to taking again a subsequence, that
the subsequence at hand is indiced by an increasing sequence $(N_p)_{p\in\NN}$ of integers verifying $j_{N_p}<i_{N_p}$ 
and such that
\bq
\sum_{p\in\NN}N_p^{-6/5}&<&+\iy\eq
Then by Borel-Cantelli Lemma and Lemma \ref{saut}, we have that a.s., $R^{(N_p,i_{N_p})}_{\tau_{N_p}}$
is diverging to infinity for large $p$. But conditionally on that, we are brought back to the 
situation where (\ref{iN5}) is satisfied (replacing $i_N$ by $R^{(N_p,i_{N_p})}_{\tau_{N_p}}$),
so Theorem \ref{th1} follows.
\begin{rem}
To avoid the above a.s.\ argument,
we can re-examine the proof of proposition~\ref{pr1},
and see that instead of deterministic initial conditions $(i_N)_{N\in\NN}$,
we could have considered random initial conditions $(i_N)_{N\in\NN}$
(such that for each $N\in\NN^*$, $i_N$ is independent from the randomness necessary to the evolution of  the Markov chain
$R^{(N,i_N)}$) diverging to infinity in probability. Thus, via lemma \ref{facile},
to obtain the wanted convergence in law, it is sufficient to show that
$R^{(N,i_N)}_{\tau_N}$ is diverging to infinity in probability.
But this is an easy consequence of lemma \ref{saut}.
\end{rem}

\section{Planar genealogical tree based representations}\label{planar}

As we promised in the end of the introduction, this last section is concerned with designing
a description of the invariant measure $\widehat{\Gamma}_{\mu}$
in terms of planar genealogical trees. Most of the arguments are only sketched, since they follow the same lines as the ones developed in~\cite{dpr}.

We start by recalling that for any pair of mappings $a\leq b$, we have $a=cb$ for some non unique $c\in\Aa$. The 
mapping $b$ induces the following equivalence relation $\sim_b$ on $\Aa$
$$
c\sim_b c'\Longleftrightarrow cb=c'b
$$
Equivalently, $c\sim_b c'$ if and only if the restriction of $c$ and $c'$ to the image of $a$ are equal, so that
the equivalence class $\overline{c}\in\Aa_{/{\sim_b}}$ can be seen as the corresponding map $\overline{c}$ from $b([N])$ into $[N]$; and we therefore have $$\left| \Aa_{/{\sim_b}}\right|=N^{|b|}$$
 In terms of backward genealogies, for a given fixed pair of mappings
 $a\leq b$,  the mapping $\overline{c}$ represents
the way the individuals $b([N])$ choose their
parents, and the composition mapping $a=\overline{c}b$ represents 
the way the individuals $[N]$ choose their
grand parents in $a([N])$. Notice that in this interpretation, not all the individuals  in the 
previous generations are ancestors, the individuals that do not leave
descendants are not counted. 

We further suppose that we are given a
function
$$
\theta~:~(a,b)\in\{(a',b')\in\Aa^2~:~a'\leq b'\}\mapsto\theta(a,b)\in\RR
$$
such that 
$
\theta(a,b)=\theta(a',b')
$, as soon as the pair of mappings $(a,a')$, and resp. $(b,b')$, only differs in the way the 
sets $a([N])$ and $a'([N])$, and  resp. $b([N])$ and $b'([N])$, are labeled. That is, $\theta (a,b)=\theta (a',b')$ whenever there exist increasing bijections $\alpha$ resp. $\beta$ from $a([N])$ to $a'([N])$ resp. $b([N])$ to $b'([N])$ s.t. $a'=\alpha a$ and $b'=\beta b$. We write $(a',b')\equiv (a,b)$ when this property is satisfied. The relation $\equiv$ is an equivalence relation.
We also mention that, for a given pair $(a,b)$, there are 
$$\left(\begin{array}{c}{N}\\{|a|}\end{array}\right)\times\left(\begin{array}{c}{N}\\{|b|}\end{array}\right)=\frac{(N)_{|a|}}{|a|!}\times
\frac{(N)_{|b|}}{|b|!}
$$ pairs of mappings s.t. $(a',b')\equiv (a,b)$.
Moreover, there is a unique pair $(a',b')$ with the property $a'([N])=[|a|]$ and $b'([N])=[|b|]$. We use the familiar combinatorial terminology and refer to this process (the replacement of the set $a([N])$ by $[|a|]$, of $a$ by $a'$, and so on), as the \it standardization process \rm (of subsets of $\NN$, of maps between these subsets...)
Therefore, to evaluate $\theta(a,b)$, up to a change of indexes we can replace
the pair of sets $(a([N]),b([N]))$, by the pair $([|a|],[|b|])$, and the pair of mappings  $(a,\overline{c})$,  by a unique pair of surjections 
$(a',c')\in\left([|a|]^{[N]}\times [|b|]^{[|a|]}\right)$.  In other terms, the pair of mappings $b\leq a$ is canonically associated to a pair of surjections 
$$
[|b|]\stackrel{c'}{\longleftarrow}[|a|]\stackrel{a'}{\longleftarrow}[N]
$$
To describe precisely the planar tree representations
of the stationary population of the neutral genetic model (\ref{defgen}), it is
convenient to introduce a series of multi index notation.
\begin{defi}
We let $\Ga_q$ be the symmetric group of all permutations
of the set $[q]$, and for any weakly decreasing sequence of integers 
$$
{\bf q}\in\Qa_n:=\left\{(q_0,\ldots,q_{n})\in [N]^{n+1}~:~q_0=N\geq q_1\geq \ldots \geq q_{n}>1\right\}
$$
we use the multi index notation
$$
{\bf q}!:=\prod_{0\leq k\leq n}q_k!\qquad (N)_{\bf q}=\prod_{0\leq k\leq n} (N)_{q_k}
\quad
\mbox{\rm and} \quad |{\bf q}|:=\sum_{0\leq k\leq n}q_k
$$
We also introduce the sets
\begin{eqnarray*}
\Ca_n({\bf q})&:=&\left\{
\cb=(c_0,\ldots,c_n)\in\Ca_n~:~\forall 0\leq p\leq n~~ |c_p|=q_{p+1}\right\}\quad\mbox{\rm and}\quad \Ga_n({\bf q}):=\prod_{0\leq p\leq n}\Ga_{q_p}
\end{eqnarray*}
with the convention $q_{n+1}=1$.  
Finally, we denote by $\Sa_n({\bf q})$
the set of all sequences of sujections  
${\bf a}=(a_0,\ldots,a_n)\in\left([q_1]^{[N]}\times\ldots\times[q_n]^{[q_{n-1}]}\times[1]^{[q_n]}\right)$. 
\end{defi}

Running back in time the  
arguments given above, 
any coalescent sequence of mappings $\cb=(c_0,\ldots,c_n)\in\Ca_n({\bf q})$ is associated in a canonical way to a unique sequence of surjective mappings 
$$
{\bf a}=(a_0,\ldots,a_n)\in
\Sa_n({\bf q})
$$
so that, up to standardization, we have that
$$
{\cb}=\pi({\bf a}):=(a_{0},a_{1}a_0,a_2a_1a_0,\ldots,a_na_{n-1}\ldots a_0)
$$
To put this again in another way, the coalescence sequence $\cb$ has, up to standardization, the following
backward representation
\begin{equation}\label{seqg}
\{1\}\stackrel{a_n}{\longleftarrow}[q_{n}]
\stackrel{a_{n-1}}{\longleftarrow}[q_{n-1}]\stackrel{}{\longleftarrow}\ldots\stackrel{}{\longleftarrow}[q_2]\stackrel{a_1}{\longleftarrow}[q_1]\stackrel{a_0}{\longleftarrow}[N]
\end{equation}
In terms of genealogical trees, the mapping $a_0$ describes the 
way the $N$ individuals in the present generation choose their
parents among $q_1$ ancestors; the mapping $a_{1}$ describes the 
way these $q_1$ individuals again choose their parents among  $q_2$ ancestors, and so on. 

From these considerations, by symmetry arguments, the invariant measure $\widehat{\Gamma}_{\mu}$ of Prop. \ref{conv} can be expressed in the following way
$$
\widehat{\Gamma}_{\mu}(F)=\sum_{n\geq 0}\sum_{{\bf q}\in\Qa_n}~\frac{1}{N^{|{\bf q}|}}~\frac{(N)_{\bf q}}{{\bf q}!}~\sum_{{\bf a}\in\Sa_n({\bf q})}\mu_{\pi(\bf a)}(F)
$$
We also have the more synthetic representation formula
$$
\widehat{\Gamma}_{\mu}(F)=
\sum_{{\bf a}\in\Sa}      \frac{1}{\rho(\bf a)!}~
\frac{(N)_{\rho({\bf a})}}{N^{|\rho({\bf a})|}}~
   \mu_{\pi(\bf a)}(F)
\quad
\mbox{\rm with}\quad 
\Sa:=\cup_{n\geq 0}\cup_{{\bf q}\in\Qa_n}\Sa_n({\bf q})
$$ 
and the multi index mapping $\rho({\bf a}):={\bf q}$, for any ${\bf a}=(a_0,\ldots,a_n)\in
\Sa_n({\bf q})$.

Although this parametrization by sequences of surjective maps of the expansion of the invariant measure is already much better than the one in Prop.~\ref{conv}, it does not take into account the full symmetry properties of $\widehat{\Gamma}_{\mu}$. To take fully advantage of it, notice first of all that the labelling of individuals in a genetic population is arbitrary. From the algebraic point of view, the invariant probability measure constructed in  Prop.~\ref{conv} inherits this property in the sense that 
$$\widehat{\Gamma}_{\mu}(F)=\widehat{\Gamma}_{\mu}(F^\sigma )$$
for any $\sigma\in \Ga_n$, where $F^\sigma (x_1,...,x_n):=F(x_{\sigma (1)},...,x_{\sigma (n)})$. In particular, when computing $\widehat{\Gamma}_{\mu}(F)$ we will always assume from now on that $F$ is symmetry invariant, that is, that $F^\sigma =F$ for all $\sigma\in \Ga_n$ -a property that we write $F\in B_b^{sym}(E^N)$.
If $F$ does not have this property, its symmetrization $\overline{F}$ has it, where $$\overline{F}:=\frac{1}{n!}\sum\limits_{\sigma\in\Ga_n}F^\sigma ,$$
and the invariance properties of $\widehat{\Gamma}_{\mu}$ insure that $\widehat{\Gamma}_{\mu}(F)=\widehat{\Gamma}_{\mu}(\overline{F})$

We then consider the following 
natural action of the permutation
group $ \Ga_{\bf q}$ on the set $\Sa_n({\bf q})$, defined for any
$ {\bf s}=(s_{p})_{0\leq p\leq n}\in \Ga_n({\bf q})$, and ${\bf a}\in \Sa_n({\bf q})$ by
$$
{\bf s}({\bf a}):=(s_{1}a_0s_0^{-1},\ldots,a_{n} s_{n}^{-1})
$$
We let  $\Za({\bf a}):=\{{\bf s}~:~{\bf s}({\bf a})={\bf a}\}$
be the stabilizer of ${\bf a}$.

This group action induces a partition of the set $\Sa_n({\bf q})$ into orbit sets or equivalent classes, where $\ab$ and $\ab'$ are equivalent
if, and only if, there exists ${\bf s}\in \Ga_n({\bf q})$ such that
$\ab'={\bf s}({\bf a})$. 
The set of equivalence classes is written $I_n({\bf q})$. In terms of ancestral lines, the genealogical trees 
associated with the sequences of mappings ${\bf s}({\bf a})$  and
${\bf a}$ only differ by a change of labels of the ancestors, at each level set. By induction on $n$, it can be shown that each orbit, or equivalence class, contains an element in
the subset $\Sa'_n({\bf q})\subset \Sa_n({\bf q})$ of all sequences of weakly increasing surjections. 
This observation can be used, together with standard techniques such as lexicographical ordering of sequences of sequences of integers to construct a canonical representative in $\Sa'_n({\bf q})$ for each equivalence class. 
Also notice the mappings
${\bf a}\mapsto \rho({\bf a})$, and ${\bf a}\mapsto \mu_{\pi({\bf a})}$ are invariant; 
that is we have that $\rho({\bf a})=\rho({\bf a}')$, and 
$\mu_{\pi({\bf a'})}(F)=\mu_{\pi({\bf a})}(F)$ for any $F\in B_b^{sym}(E^N)$, as soon as $\ab'={\bf s}({\bf a})$, for some relabeling permutation sequence ${\bf s}$. Thus, applying the class formula, we find the following formula.
\begin{prop}
For any function $F\in\Ba_b^{sym}(E^N)$, we have 
$$
\widehat{\Gamma}_{\mu}(F)=
\sum_{{\bf a}\in\Sa'}      \frac{1}{|\Za(\bf a)|}~
\frac{(N)_{\rho({\bf a})}}{N^{|\rho({\bf a})|}}~
   \mu_{\pi(\bf a)}(F)
$$
with the set $
\Sa':=\cup_{n\geq 0}\cup_{{\bf q}\in\Qa_n}I_n({\bf q})  
$. 
\end{prop}
This functional representation formula can be interpreted in terms
of genealogical trees. The precise description of this alternative interpretation
is notationally consuming, so that it will be only sketched. The reader is refered to \cite{dpr} for details on the subject.

Recall that a rooted tree ${\bf t}$ is an acyclic, connected, directed graph in which any vertex has at most an outgoing edge. 
The paths are oriented from the vertices to the root, and a leaf in a tree is a vertex without any incoming edge.
Notice first that any sequence of mappings ${\bf a}=(a_0,...,a_n)$ in $\Sa_n({\bf q})$ gives rise to a directed graph -defined as usual: that is, to any element $k$ in $[q_i]$ we associate a vertex of the graph and a directed arrow to the vertex associated to the element $a_i(k)$ in $[q_{i+1}]$. Two sequences in $S_n({\bf q})$ are equivalent under the action of $\Ga_n({\bf q})$ if and only if they have the same underlying abstract graph -see \cite{dpr}. It follows that  $I_n({\bf q})$ is 
isomorphic to the set
$\Ta_n({\bf q})$ of all rooted trees ${\bf t}$, with height $ht({\bf t})=(n+1)$, with $q_{n-(p-1)}$ vertices at each level $p=1,\ldots, (n+1)$, and  with leaves only at level $(n+1)$. 
The set of all planar trees is denoted by $ \Ta:=\cup_{n\geq 0}\cup_{{\bf q}\in\Qa_n}\Ta_n({\bf q})  
$. Next, with a slight abuse of notation, for any choice 
${\bf a}\in \Sa'_n({\bf q})$ of a representative of a tree ${\bf t}\in \Ta:=\cup_{n\geq 0}\cup_{{\bf q}\in\Qa_n}\Ta_n({\bf q})  
$ we set 
$$
\mu_{{\bf t}}=\mu_{\pi({\bf a})}\qquad \Za(\bf t)=\Za(\bf a)
\quad\mbox{\rm and}\quad
\rho({\bf t})=\rho({\bf a})
$$
In terms of genealogical trees, we finally have that
$$
\widehat{\Gamma}_{\mu}(F)=
\sum_{{\bf t}\in\Ta}      \frac{1}{|\Za(\bf t)|}~
\frac{(N)_{\rho({\bf t})}}{N^{\left|\rho({\bf t})\right|}}~
   \mu_{\bf t}(F)
   $$
A closed formula of the cardinal of the stabilizer $\Za(\bf t)$
can be easily derived using the combinatorial techniques
developed in the article~\cite{dpr}. Firstly, we recall some notions.
\begin{defi}
A forest ${\bf f}$ is a multiset of trees, that is an element of the commutative monoid 
on the set of trees. We denote by $B({\bf t})$ the forest obtained by cutting the root of tree
${\bf t}$; that is, removing its root vertex, and all its incoming edges. Conversely, we denote by $B^{-1}({\bf f})$ the tree deduced from 
the forest ${\bf f}$ by adding a common root to its rooted tree. We write 
\begin{equation}\label{normalf}
{\bf f}:={\bf t}_1^{m_1}\ldots {\bf t}_k^{m_k}
\end{equation}
for the forest with the trees ${\bf t}_i$ appearing with multiplicity
${m_i}$, with $1\leq i\leq k$.  When the trees $({\bf t}_i)_{i=1...k}$
are pairwise distinct, we say that the forest is written in normal form.
\end{defi}

\begin{defi}
The symmetry multiset of a tree is 
defined as follows
$$
{\bf t}=B^{-1}({\bf t}_1^{m_1}\ldots {\bf t}_k^{m_k})
\Longrightarrow{\bf S}({\bf t}):=(m_1,\ldots,m_k)
$$
The symmetry multiset of a forest is the disjoint union of the symmetry multisets of its trees
$$
{\bf S}({\bf t}_1^{m_1}\ldots {\bf t}_k^{m_k}):=\left(\underbrace{{\bf S}({\bf t}_1),\ldots,{\bf S}({\bf t}_1)}_{m_1-\mbox{terms}},\ldots,
\underbrace{{\bf S}({\bf t}_k),\ldots,{\bf S}({\bf t}_k)}_{m_k-\mbox{terms}}
\right)
$$
\end{defi}
Following the proof of theorem 3.8 in~\cite{dpr}, we find the following closed formula
$$
|\Za({\bf t})|=\prod_{0\leq i< ht({\bf t})}{\bf S}(B^i({\bf t}))!,
$$
where we use the multiset notation ${\bf S}({\bf t})!:=\prod\limits_{i=1}^km_i!$ if ${\bf S}({\bf t})=(m_1,...,m_k)$.
The above discussion is summarized in the following proposition.
\begin{prop}\label{mainformula}
$$
\widehat{\Gamma}_{\mu}(F)=
\sum_{{\bf t}\in\Ta}   ~   \frac{{(N)_{\rho({\bf t})}}}{
{N^{\left|\rho({\bf t})\right|}}~\prod_{0\leq i< ht({\bf t})}{\bf S}(B^i({\bf t}))!}~~
   \mu_{\bf t}(F)
   $$
\end{prop}

\end{document}